\documentclass{amsart}
\usepackage[english]{babel}
\usepackage{mathrsfs}
\usepackage{epsfig,afterpage}
\usepackage{psfrag}
\usepackage{graphicx}
\usepackage{amsmath}
\usepackage{amssymb}
\usepackage{amsfonts}
\usepackage{multirow}
\usepackage{color}

\usepackage{epigraph}

\begin{document}

\title[Remembering Igor\,R.~Shafarevich]{Remembering Igor\,R.~Shafarevich: \\ The Path to Freedom}

\author[A.O. Remizov]{A.O. Remizov}
\address{Moscow State University (Moscow, Russia)}
\email{alexey-remizov@yandex.ru}


\begin{abstract}
2023 year marks the hundredth birth anniversary of prominent Russian mathe\-ma\-ti\-cian and thinker Igor Rostislavovich Shafarevich (1923\,--\,2017).
The article presents a selection of quotations from various works of him, devoted to the future of mathematics and science in general,
with minor comments and explanations.
\end{abstract}

\maketitle

\section{Foreword}

My first (long-distance) acquaintance with Igor\,R.~Shafarevich took place more than thirty years ago, when I was a schoolboy and the magazine ``{Novy Mir}''\footnote{Russian-language monthly literary magazine since 1925. In November 1962 it became famous for publishing Aleksandr Solzhenitsyn's groundbreaking ``One Day in the Life of Ivan Denisovich''.}
published his article ``Two Roads to One Cliff'' (1989). Although the author's main idea was not actually understood by me at that time, as it did not fit with
the strong flow of information that was then pouring out in the media and in personal conversations,
the text itself struck me with beauty and logical structure and was perceived as a beautiful work of art.
If someone had said then that I was to co-author him, albeit in a math textbook\footnote{
This book \cite{LinAlg} is based on a series of lectures given in the 1950s through 1970s by I.\,R.~Shafarevich
at the Faculty of Mechanics and Mathematics of Moscow State University.},
I~would have taken it as an obvious joke somewhat similar to becoming the author of the novel ``War and Peace'' by Leo Tolstoy.

However, life turned out to be full of surprises: I had the good fortune to communicate with Shafarevich for a decade, and this communication is one of the most strong and vivid impressions of my life. Unfortunately, this vivid impression is difficult to convey on paper, so I decided to focus on one issue that I believe is of great interest.

\medskip

In some of Shafarevich's works we find reflections on mathematics and on science in general. Here we are dealing with a very rare phenomenon~-- a view of mathematics by a greatest mathematician taken from the outside.

As we know, Shafarevich was fascinated by mathematics in his school years, then externally graduated from the Faculty of Mechanics and Mathematics of the Moscow State University at the age of 17, obtained his first scientific results and became a professor six years later. I remembered the following episode: I once asked Shafarevich what I could read about an old mathematical theory of Alfred Clebsch. After a moment's thought, he named the book ``Vorlesungen \"uber Geometrie'' by Felix Klein and added that he had read it as a schoolchild and it had made a great impression on him:
{\it I had two books, which I was reading at the same time: the Klein book and ``The Three Musketeers'' by Alexandre Dumas. And I myself was surprised to find myself reading the first one all the time, and the book of Dumas was lying open on the same page} -- he said.

Certainly, all of this shows the deepest immersion in math and a degree of detachment from the realities of the world around us. A person on this path is often tempted to think of their science as of self-evident value and the world around them as just a breeding ground. However, Shafarevich did not fall into the temptation to live in an ``ivory tower''. Like Leo Tolstoy, who pondered {\it what men live by}\footnote{
The title of a short story by Leo Tolstoy written in~1881, which made a deep influence in literature and culture.},
he asks the question about the existence of ``the external meaning'' of science and offers his answer,
naturally connected with his reflections on the future of Russia, and the whole of humanity.

It is apparent to me that his thoughts are as relevant today as they were half a century ago, when they were expressed. They are exactly the ones I would like to talk about. The following quotations from several works by Shafarevich not only shed light on the reasons for the decline of science that is taking place before our eyes, but also offer, if not a way out of it, then, at least, a path on which the way out can be found.

\medskip

\section{Fateful moments}

{\begin{tabbing}
\footnotesize
Blessed is he who visits this life \\
\footnotesize
at its fateful moments of strife  \\
\end{tabbing}}
\vspace{-6mm}
{\begin{tabbing}
\emph{\footnotesize
Fyodor Tyutchev}. \\
\emph{\footnotesize
``Cicero''. }
\end{tabbing}}


\vspace{4mm}

For many years reading Shafarevich's works on various subjects, I never tired of being amazed at the diversity of issues and topics he touched upon. Gradually it became evident that all of them, to a greater or lesser extent, were directed towards one goal~-- to understand the fate of mankind and, first of all, of Russia and the Russian people. It is not accidental that one of his first works is called ``Does Russia have a future?'' (1974), and one of his last works, summarizing his many years of reflection, is ``The Future of Russia'' (2005). In these works, the author, in his own words, tries to {\it look under the lid of the pot in which our future is being boiled}.

The search led him to analyze many seemingly extraneous issues, from the Inca Empire and Antiquity to the French Revolution and the role of Protestantism in the creation of modern technological civilization. It seems to me necessary to say a few words about the background against which these reflections took place.

\medskip

After the great wars that shook the world in the middle of the 20th century, there was a fall or significant weakening of totalitarian regimes in Europe and the USSR, making the future not as fatally predetermined as it had been before. However, instead of a long-awaited rest, humanity found itself facing the threat of new global catastrophes: nuclear war, avalanche population growth and famine, destruction of nature, uncontrolled migration and national problems... Technical progress, which had served mankind (or, at least, some part of it) faithfully and loyally and seemed to be a guarantee of an imminent happy future, suddenly got out of control and began to threaten its creators.

This was noticed by many people, both in Russia and in the West, where at this time there was an explosive growth of futurological research. Under Soviet conditions, however, it was possible to discuss such issues only in the form of underground press often called ``Samizdat''\footnote{
{Samizdat} (means ``self-publishing'') was an unofficial (often censored and underground) activity across USSR and the Eastern Bloc in which individuals reproduced texts and passed them from reader to reader. The most of quotations in this article are taken from works first published in Samizdat.
}, which also grew rapidly during this period.
One of the most interesting was the work~\cite{Amal} of Andrei Amalrik\footnote{
Andrei Amalrik (1938\,--\,1980) was a Soviet writer, art collector and dissident.
He was best known for his essay ``Will the Soviet Union Survive Until 1984?'' (1970).
Soon after that, he was arrested and convicted receiving a sentence of three years in a labor camp in Kolyma.
After that, he was arrested again. In 1976, Amalrik and his wife emigrated to the Western world.
}.
The~treatise~\cite{Sahar} by Andrei Sakharov\footnote{
Andrei Sakharov (1921\,--\,1989) was a Soviet physicist and dissident. He did fundamental works in theoretical physics and contributed to the Soviet program of nuclear weapons (including the design of ``Tsar Bomba''). Also known for his political activism for individual freedom, human rights, civil liberties and reforms in USSR (Nobel Peace Prize, awarded in 1975).
} also attracted considerable attention.

At the same time, articles by Solzhenitsyn, Shafarevich and some other authors were written and published in 1974 in the collection called
``From Under the Rubble'', the initiative of creating which belonged to Solzhenitsyn.
Immediately striking is the difference between the articles in this collection and most of the texts that circulated in the under\-ground press (Samizdat), including those mentioned above. It's not just a matter of disagreement on specific issues, but it feels like a fundamentally different way of looking at life, as if we deal with Euclidean and non-Euclidean geometries. (A~direct polemic with Sakharov's treatise~\cite{Sahar} is contained in Solzhenitsyn's article~\cite{Solz} opening the collection.)

For example, estimating Amalrik's essay~\cite{Amal} as ``one of the most brilliant and clever works that Russian thought gave after the revolution'',
Shafarevich further writes \cite{Sh-Estli} the following:

\begin{quote}
{\it
The value of his book, as I see it, is that it has followed one possibility through to its logical end, that it has exhausted one train of thought. If you look at history as the product of the interaction of economic factors, or from the point of view of the interplay of the interests of different social groups and individuals, and the rights that guarantee these interests, then Russia indeed has no future~-- Amalrik's arguments are unanswerable. But there are, after all, historical processes that depend on quite other principles... Four hundred years earlier, for that matter, when an unknown monk named Luther challenged the greatest force in the world at that time, he seemed to be going counter to all social and historical laws.
}
\end{quote}

\medskip

Following Oswald Spengler, Shafarevich believed that the question of the future of Russia (as well as of any other country and, even more so, of mankind) should be considered not in the ``linear one-dimensional'' concept of progress, which prevailed in the XIX century, but on a vast ``historical field'',
on which different cultures operate independently, like living organisms,
having periods of germination, blossoming and dying. In his book ``The Decline of Western Civilization'' (written in 1918\,--\,1922), Spengler predicted the decline of Western-type civilization and its replacement by some other one\footnote{
It should be noted that some similar ideas had been expressed even earlier by the Russian scientist and thinker Nikolay Danilevsky (1822\,--\,1885) in his book ``Russia and Europe'' (1869).
}. Shafarevich discusses Spengler's concept in detail in his work~\cite{Sh-Duhovnye}.

Immediately after the publication of ``The Decline of Western Civilization'', Spengler's views were widely criticized. But now, a hundred years later, the statement about the decline of Europe can hardly be considered an unfulfilled prophecy.
However, the prediction about the decline of science seems more dubious: it was made at the beginning of the 20th century, and just a few decades later there was an explosive growth in many sciences, including mathematics (where even completely new fields appeared) and physics (suffice it to mention the discovery of the theory of relativity and quantum physics, which overturned all previous ideas about the world). And no matter how one evaluates the fruits brought by these scientific achievements in the daily lives of ordinary people, it is impossible to deny that the 20th century is the century of the rise, not the decline of science.
Shafarevich felt that this prediction of Spengler's should be seen in a slightly more distant perspective. Perhaps it is now that mankind has come close to the line that can be called the beginning of the decline of traditional science, driven by the ``Faustian spirit''~-- the irrepressible desire to know the truth. This decline is universal, despite the differences caused by the specifics of a particular country or society, see~\cite{Parshin2019} and \cite{Butler, Glaser, Nature}.

\medskip

\section{I.\,R. Shafarevich on Science (quotations)}
\smallskip

\subsection{The path to freedom begins within us (from \cite{Sh-Estli})}
$\phantom{.}$
\smallskip

At every step life presents us with choices touching upon one particular question, and that is whether to give in to force a little, to bow to pressure, or to stand our ground and straighten our backs...  Even the boldest action in these cases no longer entails imprisonment or the permanent loss of one's job. The risk is merely one of official displeasure, the loss of regular promotion and pay raises, no new television set, no bigger apartment, no official trip abroad.

A process of barter takes place in which we pay with part of our own soul that are essential to its health and survival. Our sense of self-respect and self-confidence is replaced by ruthless hostility toward others and the cunning mentality of the slave. Worst of all, life loses its aura of happiness and meaningful purpose. The price we pay is sterility in art and science, lives wasted on week-long vigils in endless queues for objects that nobody wants...

What does life offer us in exchange? As a rule the bare necessities that keep starvation at bay and enable us to feed the children. These are not in question. What is, then? I would say: values that do not lie in the material sphere. Sometimes this is quite obvious, sometimes less so. A medal, for example, neither feeds you nor keeps you warm. A large and expensive automobile soon falls apart on our roads, parking in the cities is more difficult and with speed limited by law you get to your destination no quicker than in the cheapest vehicle. A trip abroad can be important for a budding engineer's or scientist's career, but its attraction is far greater than its usefulness. An expensive new suit keeps you no warmer than an old one patched at the elbows. And so on. None of these values has a consumer significance, their meaning is quite different: they show a man's place in the hierarchy of surrounding society. Like paper money, they have no value in themselves, but are symbols of something that men value highly.

Evidently any society, in order to exist has to arrange its members in some sort of hierarchy. The hierarchy of human society reflects that society's outlook on life. The people most skilled in the activities that are highly regarded by society possess the greatest authority. Society endows such people with symbols that underline their authority~-- nose ring, gold-braided uniform or ``Chaika'' automobile\footnote{
``Chaika'' was the large and expensive automobile produced in USSR. It was commonly used to carry top ranking bureaucracy and government officials.
}.
These symbols acquire an exceptional attraction for the members of that society, persuading them to behave in the way society prefers. It is this force that is the greatest limiting factor on our present freedom. It springs not from machine guns or barbed wire, but from our own opinions, from our inward, unquestioning acceptance of the hierarchy of surrounding society, from our assumption that a high position in it really matters...

The road to freedom begins within ourselves, when we stop clawing our way up the rungs of the career ladder or of quasi-affluence. And just as we sacrifice the best part of our souls in pursuit of these will-o'-the-wisps, so when we give them up shall we find the real meaning of life.

\smallskip

\subsection{Science is turning into a race (from  \cite{Sh-Estli})}
$\phantom{.}$
\smallskip

The sheer scale, the superorganized character of modern science has been its misfortune, even its curse. There are so many scientists and their output is so great that it is impossible to read all the publications even in one narrow specialty. The scientist's horizon dwindles to a pinpoint and he exhausts himself hying to keep abreast of his countless competitors. God's design, the divine beauty of truth as revealed in science, gives way to a bundle of petty technical problems. Science becomes a race, millions speed along without the least idea of where they are going. There is still satisfaction in this race for the few with vision, who can see a few steps ahead, but the vast majority see nothing but the heels of the one immediately in front, feel nothing but the panting breath of the one treading on their heels behind.

But even if it were possible to surmount the fact that science no longer brings the satisfaction it is capable of giving, that it deforms those who practice it, there are other reasons why it cannot go on the same way indefinitely. The output of science is now doubling every ten to fifteen years, the number of scientists is growing correspondingly, and spending on science is rising at almost the same rate. This process has been going on for two hundred to two hundred fifty years, but now it is clear that it cannot go on much longer~-- for by the end of this century spending on science would exceed the whole of society's gross product. In~practice, of course, insuperable difficulties will arise long before then~-- probably in the 1980s (remember Amalrik!). In other words, development in this direction is doomed and the only question that remains is whether science can switch to another way, whereby the discovery of the truth demands neither millions of men nor billions of money, the way trodden by Archimedes and Galileo and Mendel. That is the fundamental problem now, science's life-and-death question. It will hardly be solved by those already trapped like squirrels in its treadmill. Our hopes must rest on those who have not yet been caught in its momentum.

\smallskip
\newpage

\subsection{The purpose of math (from \cite{Sh-RazvMat})}
$\phantom{.}$
\smallskip

Mathematics grows impetuously and continuously, un- disturbed by the rebuilding and crises typical of physics, constantly enriching us with new ideas and concrete facts... What is the value of the unlimited accumulation of ideas which are, in principle, equally profound? Does not mathematics become a strikingly beautiful variant of Hegel's ``idiotic infinity''?..

We might say that the development of mathematics resembles not the growth of a living organism which retains its form and alone determines its boundaries but rather the growth of a crystal or the diffusion of a gas. Both crystal and gas will spread without limit unless checked by an external obstacle. Clearly, such a development of mathematics contradicts the sense of intelligence and beauty inevitably gained from contact with the subject, just as there is an inherent contradiction in the concept of a symphony which goes on forever. Is ours the only discipline which gives rise to this problem?

I do not think that mathematics differs radically from other forms of cultural activity. True, the entities mathematics deals with are more abstract; it rejects more accidental properties. As Plato said, there is in mathematics more knowledge of pure being and fewer opinions about objects in the visible world; in mathematics ``one seems to dream of essence''. That is why tendencies dearly discernible in mathematics, while universal, are dimly discernible in other fields of knowledge. In particular, the absence of aims and form of which we spoke above is found, I think, in practically all aspects of the life of modern man...

For several centuries now man has been in the grip of feverish activity, formless and devoid of all aim and meaning save that of unlimited growth. This was called ``progress'' and for a time it became something of a substitute for religion. Its latest offspring is modern industrial society. It has been frequently pointed out that this race is self-contradictory and leads to catastrophic material consequences: to an ever-increasing tempo of life straining human capacities, to overpopulation and to the destruction of the environment. Using mathematics as an example I want to call attention to spiritual consequences which are no less destructive: all human activity loses a global aim and becomes meaningless.

The danger here is more than just negative in nature. It does not simply consist in the fact that the intense efforts of humanity, the life of its most talented members, are ultimately devoid of meaning. The full danger lies not only in our inability to predict the consequences of actions whose purpose we do not understand. Such is the spiritual constitution of mankind that it cannot for long reconcile itself to an activity whose aim and sense elude it. Here as well as in many other phenomena, what begins to operate is the mechanism of substitution: when unable to find what they need, human beings resort to a substitute...

In particular, a mathematician may seek the purpose of his work in filling the order of a state, for which he is ready to compute the trajectory of a rocket, design an eavesdropping apparatus, or, if he is exceptionally capable, plan a whole society of hybrids~-- part-man, part-computer. It is not just the souls of scientists that are mutilated by such an order of things. There appear whole areas of mathematics devoid of that divine beauty which captivates all those familiar with our discipline. More than two thousand years of history convince us that mathematics is unable to formulate the aim necessary to direct its own development. It must therefore borrow that aim from without.

Obviously, I do not intend to try to solve this profound problem which involves not only mathematics but all human endeavours. I merely wish to indicate basic directions where one might search for a solution. There are, it seems, two such directions. One could try to derive the aim of mathematics from its practical applications. However, it is hard to believe that the justification of a higher, spiritual activity could be found in a lower, material activity. In a copy of the ``Gospel of Thomas'' discovered in 1945 Jesus says with irony: ``If the flesh was made for the spirit, it is a miracle. But if the spirit was made for the body, it is a miracle of miracles.''

If, then, we reject this path, I think there remains just one possibility. The purpose of mathematics cannot be derived from an activity inferior to it but from a higher sphere of human activity, namely, religion. Clearly, it is very difficult at the present moment to see how this can happen. But it is even more difficult to conceive how mathematics can go on developing forever without knowing what it studies and why. It is bound to perish in the next generation, drowned in a flood of publications~-- but that is, after all, only the most elementary external reason.

On the other hand, the suggested solution, as history itself has proved, is in principle possible. If we again go back to the time when mathematics came into being, we see that it then knew its purpose and that it acquired that purpose by following precisely this path. Mathematics as a science came into being in the 6th century B.\,C. in the religious community of the Pythagoreans and was part of their religion. Its aim was clear. By revealing the harmony of the world as expressed in the harmony of numbers it provided a path leading towards a union with the divine. It was this lofty aim which at that time supplied the forces necessary for a scientific feat to which in principle there can be no equal. What was involved was not the discovery of a beautiful theorem, not the creation of a new branch of mathematics, but the creation of mathematics itself. Then, almost at the moment of its birth, those properties of mathematics had already come to light thanks to which general human tendencies are more clearly apparent therein than anywhere else. This is precisely the reason why at that time mathematics served as a model for the development of the fundamental principles of deductive science.

In conclusion I wish to express the hope that for this same reason mathematics can now serve as a model for the solution of the fundamental problem of our time:
{\bf To reveal a supreme religious aim and purpose for mankind's cultural activity}.

\smallskip

\subsection{Aesthetic sense (from \cite{Sh-Iz-Istorii})}
$\phantom{.}$
\smallskip

The process of fusion of science and technology is manifested not only in the influence of science on technology, but also in the opposite direction. Science employs ever more powerful and expensive technical means: grandiose accelerators that do not even fit on the territory of one country, powerful computers. A significant part of the state budget has to be spent on science~-- comparable to the cost of the army. These funds have to be systematically distributed, i.e. science becomes administered. Success in it depends on access to its technical equipment, which is in the hands of the administrators. The role of individual talent, insight is decreasing and the role of funding and organization is increasing. It changes the nature of science itself...


Certainly, ancient Greek minds could not specifically foresee the ecological crisis and other consequences of technological progress and the comprehensive application of natural-scientific ideology. But they probably felt that some of the concepts emerging then in natural science (philosophy) were by necessity taking man out of nature and pitting him against it. The call to ``conquer nature'', i.e. to perceive themselves as its enemy, would have seemed to them simply blasphemous. This may have been the motivation behind the ``Aristotelian brake'' that slowed the development of natural science for 2000~years.

And yet the heart of a physicist-mathematician can hardly accept the view of the development of, for example, physics since the 16th century (Copernicus!) as a malicious error...

Moreover, in the case of physics, this aversion has a deep internal basis. The reason is that the whole construct of ``mathematized'' physics is strikingly beautiful. It reflects, as in a mirror, the most beautiful sections of mathematics (symplectic geometry, theory of complex varieties, algebraic geometry). Such an argument may seem ``unscientific'', it appeals to the feelings.

But just in the twentieth century, humanity has made the most tragic mistakes by following general concepts (racial or class) and silencing naive immediate feelings: of pity and aversion to violence. Apparently, immediate subjective feelings are a much more reliable guide in life. Including the aesthetic sense.

\medskip

\section{Fifty years later}

Most of Shafarevich's works cited above were written half a century ago. Although this period is negligible from the point of view of fundamental historical processes, one can speak of some ``experimental testing'' of the thoughts he expressed.

Thus, recalling the {\it two directions} of possible development of mathematics~-- the path of practical applications and the path of searching for a higher goal~-- we can say that in the vast majority of cases mankind has chosen the first path, and despite some strikingly beautiful scientific achievements, in general this choice leads to the void.
On this path, not only is it impossible to find the higher purpose and meaning of humanity's cultural activity, but the ``Faustian spirit'' gradually evaporates. Scientific activity becomes similar to a commercial project, instead of searching for the truth, the scientist begins living {the sweetly exciting life of a merchant},
gradually supplanting any aesthetic sense...
Aleksei N.~Parshin\footnote{
Aleksei Parshin (1942\,--\,2022) was a Russian mathematician, specializing in algebraic geometry and arithmetic geometry.
Among his mathematical achievements one can count a generalization of local class field theory in higher dimensions and the proof of the Mordell conjecture.
}
wrote about it vividly \cite{Parshin2019}:
\begin{quote}
{\it
In recent years, while participating in the struggle against the actions of the authorities towards science in our country, I have come to realize the power of officialdom, both at home and abroad, and the ideology underlying this power... A huge society of intelligent and energetic people turned out to be absolutely helpless in the existing situation: Lilliputians defeating Gulliver.
}
\end{quote}

\medskip

A modern scientist, unless he is a rare hermit like Alexander Grothendieck or Grigori Perelman, from a young age is immersed in an atmosphere where the search for truth is replaced by constant running around, preparing applications for various grants, reports on them, writing numerous articles for these reports, actions aimed at increasing all kinds of scientific metrics, necessary not only to move up the hierarchical ladder, but also simply to survive in the scientific environment.
More and more often there is a temptation to ``consider'' speculative models that no longer reflect the features of the actual processes, to publish useless or lacking in content, repetitive articles, or even to make a direct forgery or hackwork, if it is necessary to achieve the main goal~-- success in the struggle for grants and privileges...

\medskip

Soviet writer Yuri Nagibin\footnote{
Yuri M. Nagibin (1920\,--1994) was a Soviet and Russian writer, literary critic and novelist.
}
frankly wrote in his diary notes (1951) about a similar phenomenon, except not in science, but in literature:

\begin{quote}
{\it
Just to think that inarticulate, cold, trashy scribbled pages can turn into a marvelous piece of leather on rubber, so beautifully fitted around the leg, or into a piece of the finest wool, in which one involuntarily begins to respect oneself, or some other thing... Then you stop being disgusted by inked sheets, you want to scribble many, many more. Still, I am sure that potboiler work
cause some of the most delicate and most precious brain cells to die off. Resilient when the brain is red-hot with a real holy effort,
they rot once when solving one of ``everyday tasks''.
}
\end{quote}

\medskip

Remember the words quoted from \cite{Sh-Estli}:
{\it science is turning into a race, a crowd of millions is rushing, and no one knows where to.}
One is tempted to complement these words with vivid visual image from another quotation:
each ``racer'' is driving an automobile ``Chaika'' wearing a nose ring and an embroidered uniform...
Isn't this the image of the upcoming scientific community, and perhaps of the entire ``progressive humanity''?

One can only wonder how it will all end, in what form the collapse of the current civilization will take place, and what will replace it. It is also an agonizing question that everyone decides for himself individually: in what and to what extent one can bend under the pressure of this inexorable element, what kind of a nose ring one can afford.

There can hardly be specific answers here, but for me the exemplary symbol of resistance will always be Igor Rostislavovich Shafarevich, whom I can never imagine either in an embroidered uniform or with a nose ring.

\end{document}